\documentstyle{amsppt} 
\NoRunningHeads 
\magnification=\magstep1 
\baselineskip=12pt 
\parskip=5pt 
\parindent=12pt 
\topskip=10pt 
\leftskip=0pt 
\rightskip=0pt 
\pagewidth{30pc} 
\pageheight{47pc} 
\topmatter 
\title\nofrills\bf Notes on very ample vector bundles on 
3-folds\endtitle 
\author Hidetoshi Maeda and Andrew J. Sommese\endauthor 
\abstract 
Let $\Cal E$ be a very ample vector bundle of rank two on a smooth 
complex projective threefold $X$.  An inequality about the third Segre 
class of $\Cal E$ is provided when $K_X+\det \Cal E$ is nef but not big, 
and when a suitable positive multiple of $K_X+\det \Cal E$ defines a 
morphism $X\to B$ with connected fibers onto a smooth projective curve 
$B$, where $K_X$ is the canonical bundle of $X$.  As an application, the 
case where the genus of $B$ is positive and $\Cal E$ has a global 
section 
whose zero locus is a smooth hyperelliptic curve of genus $\geq 2$ is 
investigated, and our previous result is improved for threefolds. 
\endabstract 
\endtopmatter 
\document 
\subhead\nofrills{\bf Introduction}\endsubhead \par 
In what follows, varieties are always assumed to be defined over the 
field 
$\Bbb C$ of complex numbers.\par 
Let $\Cal E$ be a very ample vector bundle of rank $n-1$ on a smooth 
projective variety $X$ of dimension $n \geq 3$.  Then there exists a 
global 
section $s\in \varGamma (\Cal E)$ whose zero locus $C = (s)_0$ is a 
smooth curve on $X$.\par 
If $K_X+\det \Cal E$ is nef and big, then there exist a birational 
morphism 
$\pi :X\to X'$ expressing $X$ as the blow-up of a smooth projective 
variety $X'$ along a finite set $B$ of points (possibly empty) and an 
ample 
vector bundle $\Cal E'$ of rank $n-1$ on $X'$ such that $\Cal E = \pi 
^*\Cal 
E'\otimes \Cal O_X(-\pi ^{-1}(B))$ and that $K_{X'}+\det \Cal E'$ is 
ample. 
The polarized pair $(X',\Cal E')$ is called the {\it first reduction} of 
$(X,\Cal E)$.  In this case, under the assumption that $C$ is a 
hyperelliptic 
curve of genus $g \geq 2$, we have proved in \cite{MS} 
$$n \geq \tau \geq \frac{(n-1)(g+t+1)+2}{2(g-1)}-1,$$ 
where $\tau$ is the nefvalue of the polarized pair $(X',K_{X'}+\det \Cal 
E')$ and $t$ is the number of exceptional divisors with respect to the 
first 
reduction morphism $\pi :X\to X'$.\par 
As a continuation of the above research, we investigate the case where 
$K_X+\det \Cal E$ is nef but not big in case $n = 3$.  According to 
\cite{BS1, Theorem 3.1}, there are five possibilities in this case.  In 
this 
article we restrict ourselves to especially the following three cases 
among them, where a suitable positive multiple of $K_X+\det \Cal E$ 
defines a morphism $X\to B$ with connected fibers onto a smooth 
projective curve 
$B$. 
\roster 
\item"(a)" $X$ is a $\Bbb P^2$-bundle over a smooth curve $B$, and $\Cal 
E_F \cong \Cal O_{\Bbb P}(2)\oplus \Cal O_{\Bbb P}(1)$ for any fiber $F$ 
of the projection $\pi :X\to B$; 
\item"(b)" $X$ is a $\Bbb P^2$-bundle over a smooth curve $B$, and $\Cal 
E_F \cong T_{\Bbb P}$ for any fiber $F$ of the projection $\pi :X\to B$, 
where $T_{\Bbb P}$ is the tangent bundle of $\Bbb P^2$; 
\item"(c)" there exists a surjective morphism $\pi :X\to B$ onto a 
smooth 
curve $B$ such that a general fiber $F$ of $\pi$ is a smooth quadric 
surface $\Bbb Q^2$ in $\Bbb P^3$ with $\Cal E_F \cong \Cal O_{\Bbb 
Q}(1)^{\oplus 2}$. 
\endroster 
We note that $K_X+\det \Cal E = \pi ^*H$ for some ample line bundle 
$H$ on $B$ in every case.  The purpose of this article is to give an 
inequality in each of the above cases concerning the third Segre class 
of 
$\Cal E$ by using the double point formula.  The precise statement of 
our 
result is as follows: 
\proclaim{Theorem 1}  Let $\Cal E$ be a very ample vector bundle of rank 
two on a smooth projective $3$-fold $X$, and assume that $(X,\Cal E)$ is 
one of {\rm (a), (b)} or {\rm (c)}.  Let $D$ denote the third Segre 
class 
$s_3(\Cal E)$ of $\Cal E$, let $q$ be the genus of $B$, and let $\rho = 
\deg 
(H-K_B)$.  Then 
\roster 
\item $D(D-20) \geq 200(q-1)+37\rho$ in case {\rm (a)}; 
\item $D(D-20) \geq 168(q-1)+21\rho$ in case {\rm (b)}; 
\item $D(D-20) \geq 168(q-1)+20\rho +k$ in case {\rm (c)}, where $k$ is 
the number of singular fibers of $\pi$. 
\endroster 
\endproclaim\noindent 
By Lemma 4, the singular fibers in case (c) are biholomorphic
to irreducible and reduced quadric surfaces in $\Bbb P^3$.\par
This article is organized as follows.  In Section 0 we collect 
preliminary material.  Section 1 is devoted to the proof of Theorem 1.   
As 
an application of Theorem 1, in Section 2 we come back to the case 
where 
$C$ is a hyperelliptic curve of genus $\geq 2$, and show that case (a) 
does 
not occur and that case (b) is very restricted when the genus of $B$ is 
positive.  This allows us to improve 
\cite{MS, Theorem 6} in case of 3-folds.  It also complements
the result \cite{LPS, Theorem 2.4} as noted in
\cite{LPS, Remark 2.5}. 
\vskip1.0cm 
\subhead\nofrills{\bf 0.  Background material}\endsubhead\par 
We use the standard notation from algebraic geometry.  The tensor 
products of line bundles are denoted additively.  The numerical 
equivalence is denoted by $\equiv$.  The pullback $i^*\Cal E$ of a 
vector 
bundle $\Cal E$ on $X$ by an embedding $i:Y\hookrightarrow X$ is denoted 
by $\Cal E_Y$.  For a vector bundle $\Cal E$ on $X$, the tautological 
line 
bundle on the projective space bundle $\Bbb P_X(\Cal E)$ associated to 
$\Cal E$ is denoted by $H(\Cal E)$.  A vector bundle $\Cal E$ on a 
projective variety $X$ is said to be {\it very ample}, if the 
tautological 
line bundle $H(\Cal E)$ on $\Bbb P_X(\Cal E)$ is very ample.  We denote 
by 
$K_X$ the canonical bundle of a smooth variety $X$.  Let $\Cal E$ be a 
vector bundle of rank $n-1$ on a smooth projective variety $X$ of 
dimension $n \geq 3$ such that there exists a global section $s\in 
\varGamma (\Cal E)$ whose zero locus $C = (s)_0$ is a smooth curve 
on $X$.  Then we should note that $K_C = (K_X+\det \Cal E)_C$.\par 
Let $X$ be a smooth projective variety of dimension $n \geq 3$, and let 
$\Cal E$ be a very ample vector bundle of rank $n-1$ on $X$.  Assume 
that 
$K_X+\det \Cal E$ is nef.  Then, by the base point free theorem, a 
suitable 
positive multiple of $K_X+\det \Cal E$ is spanned and defines a morphism 
$\pi :X\to B$ with connected fibers onto a normal projective variety 
$B$. 
Assume furthermore that $\dim B = 1$.  Then, from \cite{BS1, Theorem 
3.1} for $n = 3$ and \cite{ABW, Theorem B} for $n \geq 4$, $(X,\Cal E)$ 
is 
one of the following: 
\roster 
\item"(a)" $X$ is a $\Bbb P^{n-1}$-bundle over a smooth curve $B$, and 
$\Cal E_F\cong \Cal O_{\Bbb P}(2)\oplus \Cal O_{\Bbb P}(1)^{\oplus 
(n-2)}$ for any fiber $F$ of the projection $\pi :X\to B$; 
\item"(b)" $X$ is a $\Bbb P^{n-1}$-bundle over a smooth curve $B$, and 
$\Cal E_F\cong T_{\Bbb P}$ for any fiber $F$ of the projection $\pi : 
X\to 
B$, where $T_{\Bbb P}$ is the tangent bundle of $\Bbb P^{n-1}$; 
\item"(c)" there exists a surjective morphism $\pi :X\to B$ onto a 
smooth 
curve $B$ such that a general fiber $F$ of $\pi$ is a smooth 
hyperquadric 
$\Bbb Q^{n-1}$ in $\Bbb P^n$ with $\Cal E_F\cong \Cal O_{\Bbb 
Q}(1)^{\oplus (n-1)}$. 
\endroster 
It should be emphasized that $K_X+\det \Cal E = \pi ^*H$ for some ample 
line bundle $H$ on $B$ in every case.\par 
In all that follows, suppose that $(X,\Cal E)$ is one of (a), (b) or 
(c).   
Since 
$\Cal E$ is very ample, there exists a global section $s\in \varGamma 
(\Cal E)$ whose zero locus $C = (s)_0$ is a smooth curve on 
$X$.  Let $g$ denote the genus of $C$.  We note that the restriction 
$\pi 
_C:C\to B$ of $\pi$ to $C$ is surjective.  Hence the Riemann-Hurwitz 
formula tells us that 
$$2g-2 = d(2q-2)+r,$$ 
where $q$ is the genus of $B$, $d$ is the degree of $\pi _C$, and $r$ 
is the 
degree of the ramification divisor of $\pi _C$.  Since $d = c_{n-1} 
(\Cal E 
_F)$, we have $d = 2$ in cases (a), (c) and $d = n$ in case (b).  On 
the other 
hand, 
$$K_C = (K_X+\det \Cal E)_C = (\pi ^*H)_C = \pi _C^*H,$$ 
so that 
$$2g-2 = d\delta,$$ 
where $\delta = \deg H$.  Hence 
$$d\delta = d(2q-2)+r.$$ 
This implies that $r$ is an integer
multiple of $d$.  Let $r = d\rho$.  Then 
$$\delta = 2q-2+\rho,\tag 0.1$$ 
and we conclude that 
$$\rho = \deg (H-K_B).$$ 
It should be kept in mind that 
$$2g-2 = d(2q-2+\rho).\tag 0.2$$ 
The following lemma tells us that $\rho > 0$. 
\proclaim{Lemma 2} 
The morphism $\pi _C:C\to B$ can never be unramified. 
\endproclaim 
\demo{Proof} 
Let $X_C$ denote the fiber product $X\times _BC$ of $X$ and $C$ over 
$B$.  Then $X_C$ is connected.  Let $p_1:X_C\to X$ be the first 
projection, 
and let $p_2:X_C\to C$ be the second projection.  Note that $p_1^*\Cal 
E$ 
is ample and that $p_1^{-1}(C)$ is the zero locus of the pullback of the 
section 
$s$ defining $C$.  Thus $H^0(p_1^{-1}(C),\Bbb Z) \cong H^0(X_C,\Bbb Z) = 
\Bbb Z$ by the Lefschetz-Sommese theorem \cite{LM2, Theorem 1.1} (note 
that its proof is valid without assuming the connectedness of 
$p_1^{-1}(C)$), so that $p_1^{-1}(C)$ is connected.  On the other hand, 
$p_1^{-1}(C)$ decomposes into a curve $C'$ and a curve $C''$, where 
$C'$ is 
the image of the section of $p_2$ defined by sending $x\in C$ to $(x,x) 
\in 
X_C$.  If $\pi _C$ is unramified, then $C'\cap C'' = \emptyset$.  This 
contradicts the connectedness of $p_1^{-1}(C)$. 
\qed\enddemo 
Moreover, if $C$ is a hyperelliptic curve of genus $g \geq 2$ and if $q 
> 0$, 
then we have the following 
\proclaim{Lemma 3} 
Assume that $C$ is a hyperelliptic curve of genus $g \geq 2$.  If the 
genus 
$q$ of $B$ is positive, then $\rho$ is either $1$ or $2$. 
\endproclaim 
\demo{Proof} 
By Lemma 2 we know that $\rho > 0$.  Letting $H = K_B+\Cal H$ for some 
line bundle $\Cal H$ on $B$, we have $\rho = \deg \Cal H$ and $K_C = \pi 
_C^*(K_B+\Cal H)$.  Since $C$ is hyperelliptic, $\vert K_B+\Cal H\vert$ 
cannot take any smooth curve of positive genus as its image.  Therefore 
$\deg \Cal H \leq 2$. 
\qed\enddemo 
\vskip1.0cm 
\subhead\nofrills{\bf 1.  Proof of Theorem 1} \endsubhead\par 
From now on, throughout these notes, we assume that $n = 3$, and use the 
same notation as in Section 0.  Let 
$M = 
\Bbb P_X(\Cal E)$ be the projective space bundle associated to $\Cal E$, 
let 
$p:M\to X$ be the bundle projection, and let $H(\Cal E)$ be the 
tautological 
line bundle on 
$M$.  Then $H(\Cal E)$ is very ample.  As mentioned, we have $K_X+\det 
\Cal E = \pi ^*H$.  Thus $K_M+2H(\Cal E) = p^*(K_X+\det \Cal E) = p^*\pi 
^*H$.\par 
Let $Z$ be a general element of $\vert H(\Cal E)\vert$, let $S$ be a 
general element of $\vert H(\Cal E)_Z\vert$, and let $\varphi :S\to B$ 
be 
the restriction of $\pi \circ p$ to $S$.  Then $K_S = \varphi ^*H$.   
Let $G$ 
be a general fiber of $\pi \circ p$, and let $f = S\cap G$.  Since $S$ 
can be 
regarded as the zero  locus of a section $t\in \varGamma (M,H(\Cal 
E)^{\oplus 2})$, $f$ is the zero locus of the section $t_G\in \varGamma 
(G,H(\Cal E)_G^{\oplus 2})$.  Hence $f \ne \emptyset$, and $\varphi$ is 
surjective.  In particular, this implies that $f$ is a $1$- 
equidimensional 
smooth fiber of $\varphi$.  On the other hand, since $H^0(f,\Bbb Z) 
\cong 
H^0(G,\Bbb Z) = \Bbb Z$ by \cite{LM2, 
Theorem 1.1}, $f$ is a smooth curve in $S$ and we conclude that 
$\varphi$ has connected fibers.\par 
Now let us apply the double point formula \cite{BS2, Theorem 13.1.5} to 
$(Z,H(\Cal E)_Z)$, which tells us that 
$$e(Z)-48\chi (\Cal O_Z)+84\chi (\Cal O_S)-11K_S^2-17K_SH(\Cal 
E)_S-(K_Z+H(\Cal E)_Z)^3+D(D-20) \geq 0,$$ 
where $D = H(\Cal E)_Z^3 = H(\Cal E)^4$, and $e(Z)$ is the topological 
Euler 
characteristic of $Z$.  Then 
$$D = s_3(\Cal E) = c_1(\Cal E)^3-2c_1(\Cal 
E)c_2(\Cal E),$$ 
where $s_3(\Cal E)$ is the third Segre class of $\Cal E$.  We have 
$$K_S^2 = (\varphi ^*H)^2 = 0$$ 
and 
$$(K_Z+H(\Cal E)_Z)^3 = (K_M+2H(\Cal E))_Z^3 = (p^*\pi ^*H)_Z^3 = 0.$$ 
Let us compute $\chi (\Cal O_Z), \chi (\Cal O_S)$ and $K_SH(\Cal E)_S$. 
First, since $Z$ is the blow-up of $X$ along $C$, $Z$ is birational to 
$X$. 
Thus 
$$\chi (\Cal O_Z) = \chi (\Cal O_X) = \chi (\Cal O_B) = 1-q.$$ 
Next, since $K_S = \varphi ^*H$ and $\varphi$ has connected fibers, we 
have 
$$h^0(K_S) = h^0(H).$$ 
Furthermore, by the Kodaira vanishing theorem, 
$$0 = h^1(K_M+2H(\Cal E)) = h^1(H),$$ 
since $K_M+2H(\Cal E) = p^*\pi ^*H$.  Moreover, from the exact sequence 
$$0\to \Cal O_Z(-S)\to \Cal O_Z\to \Cal O_S\to 0,$$ 
we see that 
$$h^1(\Cal O_S) = h^1(\Cal O_Z) = q.$$ 
Therefore, the Riemann-Roch theorem applied to $(B,H)$ gives 
$$ 
\align 
2q-2+\rho = \delta = \deg H &= h^0(H)+q-1\\ 
&= h^0(K_S)+q-1 = h^0(\Cal O_S)-h^1(\Cal O_S)+h^2(\Cal O_S)+2q-2\\ 
&= \chi (\Cal O_S)+2q-2. 
\endalign$$ 
Consequently 
$$\chi (\Cal O_S) = \rho.$$ 
Finally, we have 
$$ 
\align 
K_SH(\Cal E)_S &= (\varphi ^*H)H(\Cal E)_S = (p^*\pi ^*H)_ZH(\Cal E) 
_Z^2 = 
(p^*\pi ^*H)H(\Cal E)^3\\ 
&= \delta \Bbb P_F(\Cal E_F)H(\Cal E)^3 = \delta H(\Cal E_F)^3\\ 
&= \delta (c_1(\Cal E_F)^2-c_2(\Cal E_F)) = \delta (c_1(\Cal E_F)^2-d)\\ 
&= (2q-2+\rho)(c_1(\Cal E_F)^2-d). 
\endalign$$ 
Thus, in sum, 
$$e(Z)-48(1-q)+84\rho -17(2q-2+\rho)(c_1(\Cal E_F)^2-d)+D(D-20) \geq 
0.\tag 1.0$$ 
We proceed now by cases.\par 
$(1.1)$  {\it Case} (a).  In this case we have 
$$e(Z) = e(X)+e(C) = e(\Bbb P^2)e(B)+e(C) = 3(2-2q)+2-2g.$$ 
Combining this with $(0.2)$ gives $e(Z) = 6(1-q)-d(2q-2+\rho).$  Since 
$c_1(\Cal E_F)^2 = 9$ and $d = 2$, $(1.0)$ tells us that 
$$6(1-q)-2(2q-2+\rho)-48(1-q)+84\rho -119(2q-2+\rho)+D(D-20) \geq 
0.$$ 
An easy calculation shows that 
$$D(D-20) \geq 200(q-1)+37\rho,$$ 
and the result is proved.\par 
$(1.2)$  {\it Case} (b).  Here we have $c_1(\Cal E_F)^2 = 9$ and $d = 
3$.  By 
the same argument as that in $(1.1)$, we get 
$$6(1-q)-3(2q-2+\rho)-48(1-q)+84\rho -102(2q-2+\rho)+D(D-20) \geq 
0.$$ 
Consequently 
$$D(D-20) \geq 168(q-1)+21\rho.$$\par 
$(1.3)$  {\it Case} (c).  Before proceeding with the proof, we present 
the 
following 
\proclaim{Lemma 4} 
Let $(X,\Cal E)$ be as in case {\rm (c)}.  Then every fiber of $\pi$ is 
an irreducible and reduced quadric surface in $\Bbb P^3$. 
\endproclaim 
\demo{Proof} 
We set $L = \det \Cal E$.  Then, as we pointed out, $K_X+L = \pi ^*H$. 
\par 
We first claim that every fiber $F$ of $\pi$ is irreducible and 
reduced. 
Suppose to the contrary that $F$ is not irreducible and reduced for some 
$F$.  Then we can write $F = n_1F_1+\cdots +n_rF_r$ for distinct 
integral 
surfaces $F_1, \ldots , F_r$ and positive integers $n_1, \ldots , n_r$ 
with 
$n_1+\cdots +n_r \geq 2$.  Hence a general element $T\in \vert L\vert$ 
must meet $F$ in a curve $f = n_1f_1+\cdots +n_rf_r$ for distinct 
integral curves $f_1, \ldots , f_r$.  Now we know that $Lf_i \geq 2$ for 
any $i$ from \cite{LM1, Corollary 1}.  If $Lf_i = 2$ for some $i$, then 
$f_i 
\cong \Bbb P^1$ by \cite{LM1, Corollary 2}.  Since $f_i\in \vert 
L_{F_i}\vert$, we see that $f_i\cap \operatorname{Sing}(F_i) = 
\emptyset$, where $\operatorname{Sing}(F_i)$ is the singular locus of 
$F_i$.  Thus $\operatorname{Sing}(F_i)\subset F_i-f_i$.  Since 
$\operatorname{Sing}(F_i)$ is a compact algebraic set and $F_i-f_i$ is 
affine, we conclude that $F_i$ has at most isolated singularities.  This 
implies that $F_i$ is normal, and the classification of the polarized 
surfaces of sectional genus zero applies to $(F_i,L_{F_i})$ (see for 
example \cite{BS2, Corollary 3.2.10}).  However, since $L_{F_i}E \geq 2$ 
for any rational curve $E$ on $F_i$, $(F_i,L_{F_i})$ must be $(\Bbb 
P^2,\Cal O_{\Bbb P}(2))$.  Therefore $L_{F_i}^2 = 4$.  On the other 
hand, 
$L_{F_i}^2 = L_{F_i}f_i = 2$.  This is a contradiction.  Consequently 
$Lf_i 
= L_{F_i}^2 \geq 3$ for any $i$.  We should note that $Lf = L_F^2 = 8$ 
because $\Cal E_G \cong \Cal O_{\Bbb P^1\times \Bbb P^1}(1,1)^{\oplus 
2}$ for a general fiber $G$ of $\pi$.  Moreover, since $K_X+L =\pi 
^*H$, we 
have $K_T = \pi ^*_TH$, where $\pi _T$ is the restriction of $\pi$ to 
$T$. 
Thus $T$ is a properly elliptic minimal surface, so that \cite{S, Lemma 
0.5.1} tells us that $\pi _T:T\to B$ has no multiple fibers.  As a 
direct 
result of this observation, we obtain $r = 2$ and $n_1 = n_2 = 1$.   
Since 
$K_T = \pi _T^*H$, we have $(K_T+f_i)f_i = f_i^2$ for $i = 1, 2$.   
Moreover, 
since $f = f_1+f_2$, we get $f_i(f_1+f_2) = 0$.  Thus $f_1^2 = -f_1f_2 < 
0$, since $\pi _T$ has connected fibers.  Similarly $f_2^2 < 0$.  Hence 
$f_1 \cong \Bbb P^1$ and $f_2 \cong \Bbb P^1$, and so by the same 
argument as above we have $F_1 \cong \Bbb P^2$ and $F_2 \cong \Bbb 
P^2$.  Now we know that $F = F_1+F_2$, so that $\Cal O_{F_1}(F_1+F_2) 
\cong \Cal O_{F_1}$.  Since $F_1$ meets $F_2$ in a curve, the normal 
bundle $N_{F_1/X} = \Cal O_{F_1}(F_1)$ to $F_1$ in $X$ is negative.  By 
a 
well-known theorem of Grauert (see for example \cite{BS2, Theorem 
3.2.7}) there exists a holomorphic map $p:X\to Y$ onto a normal analytic 
variety $Y$ such that $p(F_1)$ is a point, $y$, and $p$ induces a 
biholomorphism $X-F_1 \cong Y-\{y\}$.  In particular, $F_1\cap F_2$ is 
contracted.  However, this is absurd because $F_1$ has no curves with 
negative self-intersection.  Consequently every fiber $F$ of $\pi$ is 
irreducible and reduced.\par 
If $F$ is smooth, then $(F,\Cal E_F) \cong (\Bbb Q^2,\Cal O_{\Bbb 
Q}(1)^{\oplus 2})$, where $\Bbb Q^2$ is a smooth quadric surface in 
$\Bbb 
P^3$.  Let $F$ be a singular fiber of $\pi$.  We claim that $F$ is a 
singular quadric surface in $\Bbb P^3$.  To see this, take a general 
element $T\in \vert L\vert$.  Then $T$ meets $F$ in an irreducible and 
reduced curve $f$.  We note that the arithmetic genus of $f$ is one 
because $(K_T+f)f = 0$.  If $f$ is not smooth, then $f$ has a single 
singular point, so that $Lf = 1$.  This contradicts the fact that $Lf 
\geq 
2$.  Hence $f$ is smooth, and the same argument as above again shows 
that $F$ is normal.  Since $K_X+L = \pi ^*H$, we have $K_F+L_F = \Cal 
O_F$, and we conclude that $F$ is a normal Gorenstein Del Pezzo surface 
with $K_F^2 = 8$.  According to \cite{Br, Theorem 1}, $F$ is one of the 
following: 
\roster 
\item"(i)" $F$ is a singular quadric surface in $\Bbb P^3$; 
\item"(ii)" $F$ is the space obtained by blowing down the zero section 
$\Delta$ of a $\Bbb P^1$-bundle $P$ on a smooth elliptic curve $\Gamma$; 
\item"(iii)" $F$ is a rational surface with only rational double points 
as 
singularities, obtained from $\Bbb P^2$ by blowing up some number 
$\alpha \leq 8$ points (iterations allowed) then blowing down some 
number $\beta \leq \alpha$ smooth rational curves, each with 
self-intersection $-2$. 
\endroster 
In case (ii) we can write $K_P+\Cal O_P(k\Delta ) = \sigma ^*K_F$ for 
some integer $k$, where $\sigma$ is the blowing-down of $\Delta$.  Since 
$(K_P+\Cal O_P(\Delta ))_{\Delta} = K_{\Delta} = \Cal O_{\Delta}$, we 
obtain 
$$ 
\align 
\Cal O_{\Delta} = (\sigma ^*K_F)_{\Delta} &= (K_P+\Cal O_P(k\Delta 
))_{\Delta}\\ 
&= (K_P+\Cal O_P(\Delta )+\Cal O_P((k-1)\Delta ))_{\Delta} = 
\Cal O_{\Delta}((k-1)\Delta ). 
\endalign$$ 
This implies that $k = 1$.  Hence $K_P+\Cal O_P(\Delta ) = \sigma 
^*K_F$.  Set $\Delta ^2 = -r$ for some positive integer $r$.  Then $K_P 
\equiv -2\Delta -r\rho$, where $\rho$ is a fiber of the bundle 
projection 
$P\to \Gamma$.  Thus $-\sigma ^*K_F \equiv \Delta+r\rho$, so that 
$$L_F\sigma (\rho ) = -K_F\sigma (\rho ) = (-\sigma ^*K_F)\rho = 
(\Delta+r\rho )\rho = 1.$$ 
This contradicts $L_F\sigma (\rho ) \geq 2$.  In case (iii), let 
$\varphi 
:G\to \Bbb P^2$ be the composite of $\alpha$ blowing-up morphisms, and 
let $\tau :G\to F$ be the composite of $\beta$ blowing-down morphisms. 
Then $K_F^2 = K_G^2 = K_{\Bbb P}^2-\alpha = 9-\alpha$.  Hence $\alpha 
= 1$, since $K_F^2 = 8$.  Consequently $G$ is the first Hirzebruch 
surface. 
Moreover, either $\beta = 0$ or $\beta = 1$.  If $\beta = 0$, then $F = 
G$.  This is impossible because $F$ is singular.  However, the case 
$\beta 
= 1$ is also impossible because $G$ has no $(-2)$-curves.  Therefore 
case (iii) does not occur.  Consequently $F$ must be a singular quadric 
surface in $\Bbb P^3$, and the result is proved. 
\qed\enddemo\par 
We return to the proof of Theorem 1.  In case (c), let $F'$ be a 
singular fiber of $\pi$, and let $k$ denote the number of singular 
fibers of 
$\pi$.  Then, since $F'$ is a singular quadric with an isolated 
singularity 
by Lemma 4, we obtain 
$$ 
\align 
e(Z) = e(X)+e(C) &= e(X-kF')+ke(F')+e(C)\\ 
&= 4(2-2q-k)+3k+2-2g = 8(1-q)-k-d(2q-2+\rho). 
\endalign$$ 
Since $c_1(\Cal E_F)^2 = 8$ and $d = 2$, it follows from $(1.0)$ that 
$$8(1-q)-k-2(2q-2+\rho)-48(1-q)+84\rho -102(2q-2+\rho)+D(D-20) 
\geq 0.$$ 
Therefore 
$$D(D-20) \geq 168(q-1)+20\rho+k,$$ 
and we have thus proved Theorem 1. 
\vskip1.0cm 
\subhead\nofrills{\bf 2.  The case of a hyperelliptic curve} 
\endsubhead\par 
Let $\Cal E$ be a very ample vector bundle of rank two on a smooth 
projective $3$-fold $X$, and assume that $(X,\Cal E)$ is one of (a), 
(b) or 
(c).  In this section we set up the following condition: 
\roster 
\item"$(*)$" There exists a global section $s\in \varGamma (\Cal E)$ 
whose zero locus $C = (s)_0$ is a smooth hyperelliptic curve of genus $g 
\geq 2$, and the base curve $B$ has genus $q > 0$. 
\endroster 
\proclaim{Theorem 5} 
Under the assumption $(*)$, case {\rm (a)} does not occur. 
\endproclaim 
\demo{Proof} 
The proof is by contradiction.  Let $(X,\Cal E)$ be as in case (a).   
Then we 
can write $X = \Bbb P_B(\Cal V)$ for some vector bundle $\Cal V$ of rank 
three on $B$.  Since 
$\Cal E_F 
\cong 
\Cal O_{\Bbb P}(2)\oplus 
\Cal O_{\Bbb P}(1)$ for any fiber $F$ of $\pi$, $\Cal F := \pi _*(\Cal 
E\otimes(-2H(\Cal V)))$ is a line bundle on $B$, and we have an exact 
sequence 
$$0\to 2H(\Cal V)+\pi ^*\Cal F\to \Cal E\to Q\to 0$$ 
for some line bundle $Q$ on $X$.  Then $Q_F = \det \Cal E_F-2H(\Cal V)_F 
\cong 
\Cal O_{\Bbb P}(1)$.  Hence $Q = H(\Cal V)+\pi ^*L$ for some line bundle 
$L$ on $B$.  We note that $Q$ is very ample, because $\Cal E$ is very 
ample.  Let $\Cal G = \pi _*Q$.  Then, since $Q_F \cong \Cal O_{\Bbb 
P}(1)$ for any $F$, $\Cal G$ is a vector bundle of rank three on $B$ 
such 
that $(X,Q) \cong (\Bbb P_B(\Cal G),H(\Cal G))$.  Thus $\Cal G$ is 
very ample.  We can write numerically $2H(\Cal V)+\pi ^*\Cal F \equiv 
2H(\Cal G)+\mu F$ for some integer $\mu$.  Thus $\det \Cal E = 2H(\Cal 
V)+\pi ^*\Cal F+Q \equiv 3H(\Cal G)+\mu F$.  Since $K_X+\det \Cal E = 
\pi 
^*H$ for some ample line bundle $H$ on $B$, it follows from $(0.1)$ that 
$K_X+\det \Cal E \equiv (2q-2+\rho)F$.  Therefore 
$$K_X \equiv -\det \Cal E+(2q-2+\rho)F \equiv -3H(\Cal G)+(2q-2+\rho 
-\mu)F.$$ 
On the other hand, from the basic relation $K_X = -3H(\Cal G)+\pi 
^*(K_B+\det \Cal G)$, we have 
$$K_X \equiv -3H(\Cal G)+(2q-2+x)F,$$ 
where $x = \deg \Cal G$.  Consequently $\rho -\mu = x$, i\.e\., $\mu = 
\rho -x$.\par 
Now, from the above exact sequence, $c_2(\Cal E) = (2H(\Cal V)+\pi 
^*\Cal 
F)Q \equiv 2H(\Cal G)^2+\mu H(\Cal G)F$.  As observed in Section 1, $D = 
c_1(\Cal E)^3-2c_1(\Cal E)c_2(\Cal E)$.  Hence 
$$\align 
D &= (3H(\Cal G)+\mu F)^3-2(3H(\Cal G)+\mu F)(2H(\Cal G)^2+\mu H(\Cal 
G)F)\\ 
&= 15H(\Cal G)^3+17\mu H(\Cal G)^2F = 15x+17\mu\tag 2.1\\ 
&= 17\rho-2x. 
\endalign$$ 
From Theorem 1, we have $D(D-20) \geq 200(q-1)+37\rho$.  Moreover, 
under the assumption $(*)$, by Lemma 3 $\rho$ is either $1$ or $2$. 
Hence $D(D-20) \geq 37$.  This implies that $D \geq 22$.  Combining this 
with (2.1) gives $22 \leq 17\rho -2x < 17\rho$, since $x > 0$.  Thus 
$\rho = 2$.  We claim that $q = 1$.  To see this, suppose that $q \geq 
2$. 
Since $\rho = 2$, we obtain $D(D-20) \geq 274$.  Thus $D \geq 30$.  By 
using (2.1) again, we have $x = 17-(1/2)D \leq 2$.  Since $x = H(\Cal 
G)^3$ and $H(\Cal G)$ is very ample, we have a contradiction.  Therefore 
$q = 1$, and $D(D-20) \geq 74$.  Hence $D \geq 24$, and by (2.1) we get 
$x = 17-(1/2)D \leq 5$.  On the other hand, applying \cite{IT, 
Proposition 1} to $\Cal G$ gives $x \geq 7$.  This is also 
absurd. 
\qed\enddemo 
\proclaim{Theorem 6} 
Under the assumption $(*)$, let $(X,\Cal E)$ be as in case {\rm (b)}.   
Then 
$g = 4$ and $q = 1$. 
\endproclaim 
\demo{Proof} 
We can write $X = \Bbb P_B(\Cal G)$ for some vector bundle $\Cal G$ of 
rank three on $B$.  Let $x = \deg \Cal G$.  Then, since $K_X = -3H(\Cal 
G)+\pi ^*(K_B+\det \Cal G)$, we have 
$$K_X \equiv -3H(\Cal G)+(2q-2+x)F.\tag 2.2$$ 
Let $\Cal V = \pi _*(\Cal E\otimes (-H(\Cal G)))$.  Then, since $\Cal 
E_F 
\cong T_{\Bbb P}$ for any $F$, $\Cal V$ is a vector bundle of rank 
three on 
$B$, and we have an exact sequence 
$$0\to \pi ^*L\to (\pi ^*\Cal V)\otimes H(\Cal G)\to \Cal E\to 0$$ 
for some line bundle $L$ on $B$.  Let $l = \deg L$ and let $v = \deg 
\Cal V$. 
Then 
$$\det \Cal E = 3H(\Cal G)+\pi ^*\det \Cal V-\pi ^*L \equiv 3H(\Cal 
G)+(v-l)F.\tag 2.3$$ 
Thus, by (2.2), 
$$K_X+\det \Cal E \equiv (2q-2+x+v-l)F.$$ 
On the other hand, since $K_X+\det \Cal E = \pi ^*H$, it follows from 
(0.1) that 
$$K_X+\det \Cal E \equiv (2q-2+\rho)F.$$ 
Hence 
$$\rho = x+v-l.\tag 2.4$$ 
From the above exact sequence, we get 
$$c_2(\Cal E)+(\pi ^*L)\det \Cal E = c_2((\pi ^*\Cal V)\otimes H(\Cal 
G)) 
= 3H(\Cal G)^2+2(\pi ^*\det \Cal V)H(\Cal G).$$ 
Thus, by (2.3), we see that 
$$c_2(\Cal E) \equiv 3H(\Cal G)^2+2vH(\Cal G)F-lF(3H(\Cal G)+(v-l)F) = 
3H(\Cal G)^2+(2v-3l)H(\Cal G)F.$$ 
Using this and the above exact sequence again, we obtain 
$$3l = c_2(\Cal E)\pi ^*L = c_3((\pi ^*\Cal V)\otimes H(\Cal G)) = 
H(\Cal 
G)^3+(\pi ^*\det \Cal V)H(\Cal G)^2 = x+v.$$ 
Therefore by (2.4), we have $\rho = x+v-l = 3l-l = 2l$, so that $\rho$ 
is 
even.  We know from Lemma 3 that $\rho$ is either $1$ or $2$ under the 
assumption $(*)$.  Thus $\rho = 2$ and $l = 1$.  Let us compute $c_1 
(\Cal 
E)^3$ and $c_1(\Cal E)c_2(\Cal E)$.  First, by (2.3) and (2.4), 
$$ 
\align 
c_1(\Cal E)^3 = (3H(\Cal G)+(v-l)F)^3 &= 27H(\Cal G)^3+27(v-l)H(\Cal 
G)^2F = 27x+27(v-l)\\ 
&= 27(x+v-l) = 27\rho = 54. 
\endalign 
$$ 
Next, 
$$ 
\align 
c_1(\Cal E)c_2(\Cal E) &= (3H(\Cal G)+(v-l)F)(3H(\Cal 
G)^2+(2v-3l)H(\Cal G)F)\\ 
&= 9H(\Cal G)^3+(9v-12l)H(\Cal G)^2F = 9(x+v-l)-3l = 9\rho-3l = 15. 
\endalign 
$$ 
Since $D = c_1(\Cal E)^3-2c_1(\Cal E)c_2(\Cal E)$, we have $D = 24$.  On 
the other hand, Theorem 1 tells us that $D(D-20) \geq 168(q-1)+42$, 
because $\rho = 2$.  If $q \geq 2$, then $D(D-20) \geq 210$, which 
implies 
that $D \geq 28$, a contradiction.  Consequently $q = 1$, and we 
conclude 
with the aid of (0.2) that $g = 4$. 
\qed\enddemo\par 
Especially for 3-folds, Theorem 5 allows us to improve \cite{MS, Theorem 
6} as follows: 
\proclaim{Theorem 7} 
Let $\Cal E$ be a very ample vector bundle of rank two on a smooth 
projective $3$-fold $X$.  Assume that $g(X,\Cal E) = 2$.  Then $(X,\Cal 
E)$ 
is one of the following{\rm :} 
\roster 
\item $X$ is a $\Bbb P^2$-bundle over a smooth curve $B$ of genus $2$, 
and 
$\Cal E_F \cong \Cal O_{\Bbb P}(1)^{\oplus 2}$ for any 
fiber $F$ of the projection $X\to B${\rm ;} 
\item $X$ is a $\Bbb P^2$-bundle over $\Bbb P^1$, $\Cal 
E_F \cong \Cal O_{\Bbb P}(2)\oplus \Cal O_{\Bbb P}(1)$ for any fiber $F$ 
of the projection $\pi :X\to \Bbb P^1$, and $K_X+\det \Cal E$ is the 
pullback of a line bundle of degree $1$ on $\Bbb P^1$ by $\pi${\rm ;} 
\item there exists a surjective morphism $\pi :X\to B$ onto a smooth 
curve $B$ of genus $\leq 1$ such that a general fiber $F$ of $\pi$ is a 
smooth quadric surface $\Bbb Q^2$ in $\Bbb P^3$ with $\Cal E_F \cong 
\Cal 
O_{\Bbb Q}(1)^{\oplus 2}$, and $K_X+\det \Cal E$ is the pullback of a 
line 
bundle of degree $1$ on $B$ by $\pi$. 
\endroster 
\endproclaim 
 
\vskip1.0cm\newpage\flushpar 
Hidetoshi Maeda\newline 
Department of Mathematical Sciences\newline 
School of Science and Engineering\newline 
Waseda University\newline 
3-4-1 Ohkubo, Shinjuku-ku\newline 
Tokyo 169-8555\newline 
Japan\newline 
e-mail: hmaeda\@mse.waseda.ac.jp 
\vskip 1.0cm\flushpar 
Andrew J. Sommese\newline 
Department of Mathematics\newline 
University of Notre Dame\newline 
Notre Dame, INDIANA 46556\newline 
U. S. A.\newline 
e-mail: sommese\@nd.edu 

\Refs 
\widestnumber\key{ABW} 
\ref 
\key ABW 
\by M. Andreatta, E. Ballico and J. A. Wi\'sniewski 
\pages 331--340 
\paper Vector bundles and adjunction 
\yr 1992 
\vol 3 
\jour Internat\. J. Math\. 
\endref 

\ref 
\key BS1 
\by M. C. Beltrametti and A. J. Sommese 
\pages 55--74 
\paper Comparing the classical and the adjunction theoretic definition 
of scrolls 
\yr 1993 
\vol 9 
\inbook Geometry of Complex Projective Varieties, Cetraro, 
1990 
\publ Mediterranean 
\publaddr Rende 
\eds A. Lanteri, M. Palleschi and D. C. Struppa 
\bookinfo Sem\. Conf\. 
\endref 

\ref 
\key BS2 
\by M. C. Beltrametti and A. J. Sommese 
\yr 1995 
\vol 16 
\book The Adjunction Theory of Complex Projective Varieties 
\publ de Gruyter 
\publaddr Berlin 
\bookinfo de Gruyter Exp\. Math\. 
\endref 

\ref 
\key Br 
\by L. Brenton 
\pages 117--124 
\paper On singular complex surfaces with negative canonical bundle, with 
applications to singular compactifications of $\Bbb C^2$ and to 
$3$-dimensional rational singularities 
\yr 1980 
\vol 248 
\jour Math\. Ann\. 
\endref 

\ref 
\key IT 
\by P. Ionescu and M. Toma 
\pages 633--643 
\paper On very ample vector bundles on curves 
\yr 1997 
\vol 8 
\jour Internat\. J\. Math\. 
\endref 

\ref 
\key LM1 
\by A. Lanteri and H. Maeda 
\pages 181--199 
\paper Adjoint bundles of ample and spanned vector bundles on algebraic 
surfaces 
\yr 1992 
\vol 433 
\jour J. Reine Angew\. Math\. 
\endref 

\ref 
\key LM2 
\by A. Lanteri and H. Maeda 
\pages 247--259 
\paper Ample vector bundle characterizations of projective bundles and 
quadric fibrations over curves 
\yr 1996 
\inbook Higher Dimensional Complex Varieties, Trento, 1994 
\publ Walter de Gruyter 
\publaddr Berlin-New York 
\eds M. Andreatta and T. Peternell 
\endref 

\ref 
\key LPS 
\by A. Lanteri, M. Palleschi and A. J. Sommese 
\pages 2307--2324 
\paper On the adjunction mapping of very ample vector
bundles of corank one 
\yr 2004 
\vol 356 
\jour Trans\. Amer\. Math\. Soc\. 
\endref 

\ref 
\key MS 
\by H. Maeda and A. J. Sommese 
\pages 74--80 
\paper Very ample vector bundles of curve genus two 
\yr 2002 
\vol 79 
\jour Arch\. Math\. (Basel) 
\endref 

\ref 
\key S 
\by A. J. Sommese 
\pages 211--220 
\paper On the nonemptiness of the adjoint linear system of a hyperplane 
section of a threefold 
\yr 1989 
\vol 402 
\jour J. Reine Angew\. Math\. 
\endref 
\endRefs
\enddocument